\documentclass[english,12pt,a4paper,twoside]{smfart}
\usepackage{amsfonts,smfthm}
\usepackage{amsmath}
\usepackage{stmaryrd}
\usepackage{amssymb}
\usepackage{graphics}
\usepackage{babel}
\usepackage{enumerate}
\usepackage{euscript}

\usepackage{newlfont}
\usepackage{latexsym}
\usepackage{graphicx}
\usepackage{a4wide}

\usepackage{cite}        
\usepackage{url}         

\title{Inverse zero-sum problems \\ and algebraic invariants}
\author{Benjamin Girard}

\address{Centre de
Math\'{e}matiques Laurent Schwartz, UMR $7640$ du CNRS, \'Ecole
polytechnique, $91128$ Palaiseau cedex, France.}
\email{benjamin.girard@math.polytechnique.fr}

\theoremstyle{plain}
\newtheorem{theorem}{Theorem}[section]
\newtheorem{proposition}[theorem]{Proposition}
\newtheorem{lem}[theorem]{Lemma}

\newtheorem{conjecture}[theorem]{Conjecture}

\theoremstyle{definition}

\theoremstyle{remark}

\def\cnp[#1,#2]{\begin{pmatrix} #1 \\#2 \end{pmatrix}}
\begin{document}

\maketitle 

\begin{abstract}
In this article, we study the maximal cross number of long zero-sumfree sequences in a finite Abelian group. Regarding this inverse-type problem, we formulate a general conjecture and prove, among other results, that this conjecture holds true for finite cyclic groups, finite Abelian $p$-groups and for finite Abelian groups of rank two. Also, the results obtained here enable us to improve, via the resolution of a linear integer program, a result of W.~Gao and A.~Geroldinger concerning the minimal number of elements with maximal order in a long zero-sumfree sequence of a finite Abelian group of rank two.  
\end{abstract}

\vspace{-0.7cm}
\section{Introduction}
\label{Introduction}
Let $G$ be a finite Abelian group, written additively. By $\exp(G)$ we denote the exponent of $G$. If $G$ is cyclic of order $n$, it will be denoted by $C_n$. In the general case, we can decompose $G$ (see for instance \cite{Samuel}) as a direct product of cyclic groups $C_{n_1} \oplus \dots \oplus C_{n_r}$ where $1 < n_1 \text{ } |\text{ } \dots \text{ }|\text{ } n_r \in \mathbb{N}$. 

\medskip
In this paper, any finite sequence $S=\left(g_1,\dots,g_\ell\right)$ of $\ell$ elements from $G$ will be called a \em sequence \em in $G$ with \em length \em $\left|S\right|=\ell$. Given a sequence $S=\left(g_1,\dots,g_\ell\right)$ in $G$, we say that
$s \in G$ is a \em subsum \em of $S$ when it lies in the following set,
$$\Sigma(S)=\left\{ \displaystyle\sum_{i \in I} g_i \text{ }|\text{ } \emptyset \varsubsetneq I
\subseteq \{1,\dots,\ell\}\right\}.$$

\noindent If $0$ is not a subsum of $S$, we say that $S$ is a \em
zero-sumfree sequence\em. If $\sum^\ell_{i=1} g_i=0,$ then $S$ is said
to be a \em zero-sum sequence\em. If moreover one has $\sum_{i \in
I} g_i \neq 0$ for all proper subsets $\emptyset \subsetneq I
\subsetneq \{1,\dots,\ell\}$, $S$ is called a \em minimal zero-sum
sequence\em.

\medskip
In a finite Abelian group $G$, the order of an element $g$ will be
written $\text{ord}(g)$ and for every divisor $d$ of the exponent of
$G$, we denote by $G_d$ the subgroup of $G$ consisting of all
elements of order dividing $d$: $$G_d=\left\{x \in G \text{ } |
\text{ } dx=0\right\}.$$ 

\medskip
\noindent For every divisor $d$ of $\exp(G)$, and every sequence $S$ in $G$, we denote by $\alpha_d$ the number of elements, counted with multiplicity, contained in $S$ and of order $d$. Although the quantity $\alpha_d$ clearly depends on $S$, we will not emphasize this dependence in the present paper, since there will be no risk of confusion.

\medskip
Let $\mathcal{P}$ be the set of prime numbers. Given a positive integer $n \in
\mathbb{N}^*=\mathbb{N}\backslash\{0\}$, we denote by
$\mathcal{D}_n$ the set of its positive divisors and we set 
$\tau(n)=\left|\mathcal{D}_n\right|$. If $n > 1$, we denote by 
$P^-(n)$ the smallest prime element of $\mathcal{D}_n$, and we put 
by convention $P^-(1)=1$. For every prime $p \in \mathcal{P}$, $\nu_p(n)$
will denote the $p$-adic valuation of $n$.

\medskip 
Let $G \simeq C_{n_1} \oplus \dots \oplus C_{n_r},$ with $1 < n_1
\text{ } |\text{ } \dots \text{ }|\text{ } n_r \in \mathbb{N}$, be a finite Abelian 
group. We set:
$$\mathsf{D}^{*}(G)=\displaystyle\sum^r_{i=1} (n_i-1) + 1 \hspace{0.2cm} \text{ as well as } \hspace{0.2cm} \mathsf{d}^{*}(G)=\mathsf{D}^{*}(G)-1.$$
By $\mathsf{D}(G)$ we denote the smallest integer $t \in \mathbb{N}^*$ such 
that every sequence $S$ in $G$ with $|S| \geq t$ contains a non-empty zero-sum subsequence.
The number $\mathsf{D}(G)$ is called the \em Davenport constant
\em of the group $G$. 

\medskip
\noindent By $\mathsf{d}(G)$ we denote the largest integer $t \in \mathbb{N}^*$ such that there exists a zero-sumfree sequence $S$ in $G$ with $|S|=t$. It can be readily seen that, for every finite Abelian group $G$, one has $\mathsf{d}(G)=\mathsf{D}(G)-1$. 

\medskip
If $G \simeq C_{\nu_1} \oplus \dots \oplus
C_{\nu_s}$, with $\nu_i > 1$ for all $i \in \llbracket 1,s
\rrbracket$, is the longest possible decomposition of $G$ into a
direct product of cyclic groups, then we set:
$$\mathsf{k}^{*}(G)=\displaystyle\sum_{i=1}^s
\frac{\nu_i-1}{\nu_i}.$$ 

\medskip
\noindent The \em cross number \em of a sequence
$S=(g_1,\dots,g_\ell)$, denoted by $\mathsf{k}(S)$, is then defined by:
$$\mathsf{k}(S)=\displaystyle\sum_{i=1}^{\ell} \frac{1}{\text{ord}(g_i)}.$$
The notion of cross number was introduced by U.~Krause in \cite{Krause84} (see also \cite{Krause91}). 
Finally, we define the so-called \em little cross number \em
$\mathsf{k}(G)$ of $G$:
$$\mathsf{k}(G)=\max\{\mathsf{k}(S) | S \text{ zero-sumfree sequence of } G \}.$$

\medskip 
Given a finite Abelian group $G$, two elementary constructions
(see \cite{GeroKoch05}, Proposition $5.1.8$) give
the following lower bounds:
$$\mathsf{D}^{*}(G) \leq \mathsf{D}(G) \hspace{0.5cm} \text{ and } \hspace{0.5cm} \mathsf{k}^{*}(G) \leq \mathsf{k}(G).$$

\medskip
The invariants $\mathsf{D}(G)$ and $\mathsf{k}(G)$ play a key r\^{o}le in the theory of non-unique factorization 
(see for instance Chapter $9$ in \cite{Narkie04}, the book \cite{GeroKoch05} which presents the different aspects of the theory, and the survey \cite{GeroKoch06} also). They have been extensively studied during last decades and even if numerous results were proved (see Chapter $5$ of the book \cite{GeroKoch05}, \cite{GaoGero06} for a survey with many references on the subject, and \cite{Girard08} for recent results on the cross number of finite Abelian groups), their exact values are known for very special types of groups only. In the sequel, we will need some of these values in the case of finite Abelian $p$-groups and finite Abelian groups of rank two, so we gather them in the following theorem (see \cite{GeroCross94}, \cite{Olso69a} and \cite{Olso69b}). 

\newpage
\begin{theorem} \label{Proprietes D et Cross number} The following two statements hold.
\begin{itemize}
\item[$(i)$] Let $p \in \mathcal{P}$, $r \in \mathbb{N}^*$ and
$ a_1 \leq \dots \leq a_r$, where $a_i \in
\mathbb{N}^*$ for all $i \in \llbracket 1,r \rrbracket$. Then, for
the $p$-group $G \simeq C_{p^{a_1}} \oplus \cdots \oplus
C_{p^{a_r}}$, we have:
$$\mathsf{D}(G)=\displaystyle\sum^r_{i=1} \left(p^{a_i}-1\right)+1=\mathsf{D}^{*}(G) \hspace{0.23cm} \text{ and }
\hspace{0.31cm} \mathsf{k}(G)=\displaystyle\sum^r_{i=1} \left(\frac{p^{a_i}-1}{p^{a_i}}\right)=\mathsf{k}^{*}(G).$$

\item[$(ii)$] For every $m,n \in \mathbb{N}^*$, we have: 
$$\mathsf{D}(C_m \oplus C_{mn})=m+mn-1=\mathsf{D}^*(C_m \oplus C_{mn}).$$ 
In particular, we have $\mathsf{D}(C_n)=n.$
\end{itemize}
\end{theorem}

\medskip
The aim of this paper is to study some inverse zero-sum problems of a special type. Instead of trying to characterize explicitly, given a finite Abelian group, the structure of long zero-sumfree sequences (see \cite{GaoGero02}, \cite{ChenSav07}, \cite{GaoGeroSchmid07}, \cite{Wolfgang} and \cite{GaoGeroGrynkie08}), or the structure of zero-sumfree sequences with large cross number (see \cite{GeroSchnCross97}), we study to what extent a zero-sumfree sequence can be extremal in both directions simultaneously. For instance, what is the maximal cross number of a long zero-sumfree sequence? Regarding this problem, we propose the following general conjecture. 

\begin{conjecture}\label{Conjecture Cross number} Let $G \simeq C_{n_1} \oplus \dots \oplus C_{n_r},$ with $1 < n_1
\text{ } |\text{ } \dots \text{ }|\text{ } n_r \in \mathbb{N}$, be a finite Abelian group.
Given a zero-sumfree sequence $S$ in $G$ verifying $\left| S \right| \geq \mathsf{d}^*(G)$, one always 
has the following inequality:
\begin{eqnarray*} \mathsf{k}(S) \leq \displaystyle\sum_{i=1}^r \left(\frac{n_i-1}{n_i}\right). \end{eqnarray*}
In particular, one has $\mathsf{k}(S) < r$.
\end{conjecture}

\medskip
One can notice that Conjecture \ref{Conjecture Cross number} is closely related to the distribution of the orders of elements in a long zero-sumfree sequence. As we will see in this paper, it provides, when it holds, useful informations on this question. In the following proposition, we gather what is currently known, to the best of our knowledge, on the structure of long zero-sumfree sequences in finite Abelian groups of rank two. This result, due to W.~Gao and A.~Geroldinger, can be found under a slightly different form in \cite{GeroKoch05}, Proposition $5.8.4$. 

\begin{proposition}\label{ordreinseqmini} Let $G \simeq C_m \oplus C_{mn}$, where $m,n \in \mathbb{N}^*$, be a finite Abelian group of rank two. For every zero-sumfree sequence $S$ in $G$ with $\left|S\right|=\mathsf{d}(G)=m+mn-2$, the following two statements hold.

\begin{itemize}
\item[$(i)$] For every element $g \in S$, one has $m \text{ } | \text{ } \text{\em ord\em}(g) \text{ } | \text{ } mn$.

\item[$(ii)$] The sequence $S$ contains at least 
$$m+mn-n\left(\frac{2m-2}{P^-(n)}+1\right)-1 \geq m-1$$ elements with order $mn$. 
\end{itemize}
\end{proposition}

The problem of the exact structure of a long zero-sumfree sequence in groups of the form $G \simeq C_m \oplus C_{mn}$ is also closely related to an important conjecture in additive group theory, which bears upon the so-called Property B. Let $n \geq 2$ be an integer. We say that $n$ has Property B if every zero-sumfree sequence in $G \simeq C_n \oplus C_n$ with $|S|=\mathsf{d}(G)=2n-2$ contains some element repeated at least $n-2$ times. 

\medskip
Property B was introduced and first studied in \cite{GaoGero99} (see also \cite{GeroKoch05}, Section $5.8$, \cite{Lettl07} 
and \cite{GaoGeroSchmid07}). It is conjectured that every integer $n \geq 2$ has Property B, and recently, it was proved that an integer $n \geq 2$ has Property B if each prime divisor of $n$ has this property (see \cite{GaoGero03}, Section $8$ and \cite{GaoGeroGrynkie08}). Therefore, it remains to solve this problem for prime numbers. Regarding this, it can be shown that Property B holds for $n=2,3,5,7$ (see \cite{GaoGero03}, Proposition $4.2$), for $n=11,13,17,19$ (see \cite{BhoHalSch08}), and consequently for every integer $n$ being representable as a product of these numbers.

\medskip
Moreover, W.~Schmid proved in \cite{Wolfgang} that if some integer $m \geq 2$ has Property B, then the zero-sumfree sequences in $G \simeq C_m \oplus C_{mn}$ with length $\mathsf{d}(G)=m+mn-2$ can be characterized explicitly for all $n \in \mathbb{N}^*$. This result provides a unified way to prove Theorem $3.3$ in \cite{GaoGero02} and Theorem in \cite{ChenSav07}. It also implies, assuming that Property B holds for every integer $n \geq 2$, that Conjecture \ref{Conjecture Cross number} holds true for every finite Abelian group of rank two.

\section{New results and plan of the paper}
\label{Section new results}
In this article, we prove that Conjecture \ref{Conjecture Cross number} holds for several types of finite Abelian groups. To begin with, in Section \ref{Section premiers resultats}, we prove some consequences of this conjecture in the cases where it holds. For instance, Conjecture \ref{Conjecture Cross number}, if true, would imply simultaneously two classical and long-standing conjectures related to the Davenport constant of finite Abelian groups of the form $C^r_n$.
\begin{proposition}\label{La Conjecture implique la Conjecture sur Davenport} Let $n,r \in \mathbb{N}^*$ be such that Conjecture \ref{Conjecture Cross number} holds for the group $C^r_n$. Then, one has the following equality:
$$\mathsf{D}(C^r_n)= r(n-1)+1.$$
Moreover, every zero-sumfree sequence $S$ in $C^r_n$ with $\left|S\right|=\mathsf{d}(C^r_n)=r(n-1)$ consists only of elements with order $n$. 
\end{proposition}

More generally, Conjecture \ref{Conjecture Cross number}, if true, would provide the following general upper bound for the Davenport constant of a finite Abelian group.
\begin{proposition}\label{La Conjecture donne une borne pour la constante de Davenport} Suppose that Conjecture \ref{Conjecture Cross number} holds for $G \simeq C_{n_1} \oplus \dots \oplus C_{n_r},$ with $1 < n_1 \text{ } |\text{ } \dots \text{ }|\text{ } n_r \in \mathbb{N}$. Then, one has the following inequality:
$$\mathsf{D}(G) \leq \displaystyle\sum_{i=1}^r \frac{n_r}{n_i}(n_i-1)+1=\mathsf{D}^*(G)+\displaystyle\sum_{i=1}^r \left(\frac{n_r}{n_i}-1\right)\left(n_i-1\right).$$
\end{proposition}

Then, in Section \ref{Section premiers resultats} also, we prove that Conjecture \ref{Conjecture Cross number} holds true for finite cyclic groups and finite Abelian $p$-groups.
\begin{proposition}\label{La conjecture est verifiee pour les groups cycliques et les p-groupes} 
Conjecture \ref{Conjecture Cross number} holds for the following groups $G$.
\begin{itemize}
\item[$(i)$] $G$ is a finite cyclic group.

\item[$(ii)$] $G$ is a finite Abelian $p$-group.
\end{itemize}
\end{proposition}

\medskip 
In Section \ref{Outline of the method}, we present a general method which was introduced in \cite{Girard08} so as to study the cross number of finite Abelian groups. Then, using this method, we prove in Section \ref{Section lemmes} two important lemmas, which will be useful in the study of the special case of finite Abelian groups of rank two.   

\medskip
In Section \ref{Demonstrations}, we prove the two main theorems of this paper. The first one states that Conjecture \ref{Conjecture Cross number} holds for every finite Abelian group of rank two. As already mentioned in Section \ref{Introduction}, this result supports Property B (see \cite{Wolfgang}). 

\begin{theorem}\label{main theorem} Let $G \simeq C_{m} \oplus C_{mn},$ where $m,n \in \mathbb{N}^*$, be a finite Abelian group of rank two. For every zero-sumfree sequence $S$ in $G$ with $\left| S \right| \geq \mathsf{d}^*(G)= m+mn-2$, the following inequality holds:
\begin{eqnarray*} \mathsf{k}(S) \leq \left(\frac{m-1}{m}\right)+\left(\frac{mn-1}{mn}\right). \end{eqnarray*}
In particular, one always has $\mathsf{k}(S) < 2$.
\end{theorem}

The second theorem, which is proved in Section \ref{Demonstrations} as well, is an effective result which states that, in a finite Abelian group of rank two, most of the elements of a long zero-sumfree sequence must have maximal order. This result improves significantly the statement of Proposition \ref{ordreinseqmini} $(ii)$. 

\begin{theorem}\label{main theorem bis} 
Let $G \simeq C_m \oplus C_{mn}$, where $m,n \in \mathbb{N}^*$, be a finite Abelian group of rank two. For every zero-sumfree sequence $S$ in $G$ with $\left|S\right|=\mathsf{d}(G)=m+mn-2$, the following two statements hold.
\begin{itemize}
\item[$(i)$] If $n$ is a prime power, then $S$ contains at least $mn-1$ elements with order $mn$.

\item[$(ii)$] If $n$ is not a prime power, then $S$ contains at least 
$$\left\lceil \frac{4}{5}mn + \frac{(n-5)}{5}\right\rceil$$ elements with order $mn$. 
\end{itemize}
\end{theorem}

It may be observed that for every group $G \simeq C_m \oplus C_{mn}$, where $m,n \in \mathbb{N}^*$ and $n \geq 2$, there exists a zero-sumfree sequence $S$ in $G$ with $\left|S\right|=\mathsf{d}(G)=m+mn-2$, and which does not contain strictly more than $mn-1$ elements with order $mn$. Indeed, let $(e_1, e_2)$ be a basis of $G$, with $\text{ord}(e_1)=m$ and $\text{ord}(e_2)=mn$. Then, it suffices to consider the zero-sumfree sequence $S$ consisting of the element $e_1$ repeated $m-1$ times and the element $e_2$ repeated $mn-1$ times. From this point of view, Theorem \ref{main theorem bis} proves to be "nearly optimal". In addition, the general method presented in Section \ref{Outline of the method} can be successfully used to prove an analogue of Theorem \ref{main theorem bis} in the case of finite Abelian $p$-groups (see \cite{Girard09bis}).

\medskip
Finally, in Section \ref{concluding remark}, we will present and discuss a general conjecture concerning the maximal possible length of a zero-sumfree sequence with large cross number, which can be seen as a dual version of Conjecture \ref{Conjecture Cross number}.

\bigskip

\section{Proofs of Propositions \ref{La Conjecture implique la Conjecture sur Davenport}, \ref{La Conjecture donne une borne pour la constante de Davenport} and \ref{La conjecture est verifiee pour les groups cycliques et les p-groupes}} \label{Section premiers resultats}
To start with, we prove the two corollaries of Conjecture \ref{Conjecture Cross number} announced in Section \ref{Section new results}.
\begin{proof}[Proof of Proposition \ref{La Conjecture implique la Conjecture sur Davenport}] Let $S$ be a zero-sumfree sequence in $G \simeq C^r_n$ with maximal length $\left| S \right|=\mathsf{d}(G)=\mathsf{D}(G)-1$. Then, one has the following inequality: 
$$\frac{\mathsf{D}(G)-1}{n} = \frac{\left|S\right|}{n}  \leq \mathsf{k}(S) \leq r\left(\frac{n-1}{n}\right),$$
which implies that $\mathsf{D}(G) \leq r(n-1)+1=\mathsf{D}^{*}(G)$, and since $\mathsf{D}^*(G) \leq \mathsf{D}(G)$ always holds, the equality follows.
Consequently, one has: 
$$\mathsf{k}(S)= r\left(\frac{n-1}{n}\right) = \frac{\mathsf{D}(G)-1}{n} = \frac{\left|S\right|}{n},$$
and so, every element $g$ of $S$ verifies $\text{ord}(g)=\exp(G)=n$. 
\end{proof}

\begin{proof}[Proof of Proposition \ref{La Conjecture donne une borne pour la constante de Davenport}]
Let $S$ be a zero-sumfree sequence in $G \simeq C_{n_1} \oplus \dots \oplus C_{n_r},$ with $1 < n_1 \text{ } |\text{ } \dots \text{ }|\text{ } n_r \in \mathbb{N}$, such that $\left| S \right|=\mathsf{d}(G)=\mathsf{D}(G)-1$. Then, one has the following inequality: 
$$\frac{\mathsf{D}(G)-1}{n_r} = \frac{\left|S\right|}{n_r}  \leq \mathsf{k}(S) \leq \displaystyle\sum^r_{i=1} \left(\frac{n_i-1}{n_i}\right),$$
which implies the desired result.
\end{proof}

We prove now that Conjecture \ref{Conjecture Cross number} holds true for finite cyclic groups and finite Abelian $p$-groups. 
\begin{proof}[Proof of Proposition \ref{La conjecture est verifiee pour les groups cycliques et les p-groupes}]
\begin{itemize}
\item[$(i)$] Let $n \geq 2$ be an integer and let $S$ be a zero-sumfree sequence in $C_n$ with $\left|S\right| \geq \mathsf{d}^*(C_n)=n-1$. Then, it is well-known (see for instance \cite{GeroKoch05}, Theorem $5.1.10$ $(i)$) that there exists $g \in C_n$ with $\text{ord}(g)=n$ such that $S$ is of the following form: 
$$S=(\underbrace{g,\dots,g}_{n-1 \text{ times}}).$$ 
Consequently, we obtain:
$$\mathsf{k}(S)=\frac{n-1}{n},$$
which gives the desired result.

\medskip
\item[$(ii)$] Let $p \in \mathcal{P}$, $r \in \mathbb{N}^*$, and $G \simeq C_{p^{a_1}} \oplus \cdots \oplus
C_{p^{a_r}}$, with $a_1 \leq \dots \leq a_r$ and $a_i \in \mathbb{N}^*$ for all $i \in \llbracket 1,r \rrbracket$, be a $p$-group. By Theorem \ref{Proprietes D et Cross number} $(i)$, one has: 
$$\mathsf{k}(G)=\displaystyle\sum_{i=1}^r \left(\frac{p^{a_i}-1}{p^{a_i}}\right)=\mathsf{k}^*(G).$$
Then, for every zero-sumfree sequence $S$ in $G$, in particular for those verifying $\left| S \right| \geq \mathsf{d}^*(G)$, one indeed has, by the very definition of the little cross number:
$$\mathsf{k}(S) \leq \mathsf{k}(G)=\displaystyle\sum_{i=1}^r \left(\frac{p^{a_i}-1}{p^{a_i}}\right),$$
and the proof is complete.
\end{itemize}
\end{proof}

\section{Outline of a new method}
\label{Outline of the method}
Let $G$ be a finite Abelian group, and let $S$ be a sequence of elements in $G$. The general method that we will use in this paper (see also \cite{Girard08} and \cite{Girard09bis} for applications of this method in two other contexts), consists in considering, for every $d',d \in \mathbb{N}$ such that $1 \leq d' \text{ } | \text{ } d \text{ } | \text{ } \exp(G)$, the following exact sequence:
$$0 \rightarrow G_{d/d'} \hookrightarrow G_d \overset{\pi_{(d',d)}}\rightarrow \frac{G_d}{G_{d/d'}} \rightarrow 0.$$
Now, let $U$ be the subsequence of $S$ consisting of all the elements whose order divides $d$. If, for some $1 \leq d' \text{ } | \text{ } d \text{ } | \text{ } \exp(G)$, it is possible to find sufficiently many disjoint non-empty zero-sum subsequences in $\pi_{(d',d)}(U)$, that is to say sufficiently many disjoint subsequences in $U$ the sum of which are elements of order dividing $d/d'$, then $S$ cannot be a zero-sumfree sequence in $G$.

\medskip
So as to make this idea more precise, we proposed in \cite{Girard08} to introduce the following number, which can be seen as an extension of the classical Davenport constant. 

\medskip
Let $G \simeq C_{n_1} \oplus \dots \oplus C_{n_r},$ with $1 < n_1 \text{ } |\text{ } \dots \text{ }|\text{ } n_r \in \mathbb{N}$, be a finite Abelian group and $d',d \in \mathbb{N}$ be two integers such that $1 \leq d' \text{ } | \text{ } d \text{ } | \text{ } \exp(G)$. By $\mathsf{D}_{(d',d)}(G)$ we denote the smallest integer $t \in
\mathbb{N}^*$ such that every sequence $S$ in $G_d$ with $|S| \geq t$ contains a subsequence with sum in $G_{d/d'}$.

\medskip
Using this definition, we can prove the following simple lemma, which is one possible illustration of the idea we presented. This result will be useful in Section \ref{Section lemmes} and states that given a finite Abelian group $G$, there exist strong constraints on the way the orders of elements have to be distributed within a zero-sumfree sequence. 

\begin{lem}\label{lemme clef} Let $G$ be a finite Abelian group and $d',d \in \mathbb{N}$ be two integers such that $1 \leq d' \text{ } | \text{ } d \text{ } | \text{ } \exp(G)$. Given a sequence $S$ of elements in $G$, we will write $T$ for the subsequence of $S$ consisting of all the elements whose order divides $d/d'$, and we will write $U$ for the subsequence of $S$ consisting of all the elements whose order divides $d$ (In particular, one has $T \subseteq U$).  
Then, the following condition implies that $S$ cannot be a zero-sumfree sequence:
$$\left|T\right|+\left\lfloor\frac{\left|U\right|-\left|T\right|}{\mathsf{D}_{(d',d)}(G)}\right\rfloor \geq \mathsf{D}_{\left(\frac{d}{d'},\frac{d}{d'}\right)}(G).$$
\end{lem}

\begin{proof} Let us set $\mathsf{\Delta}=\mathsf{D}_{\left(\frac{d}{d'},\frac{d}{d'}\right)}(G)$. When it holds, this inequality implies that there are $\mathsf{\Delta}$ disjoint subsequences $S_1,\dots,S_{\mathsf{\Delta}}$ of $S$, the sum of which are elements of order dividing $d/d'$. Now, by the very definition of $\mathsf{D}_{\left(\frac{d}{d'},\frac{d}{d'}\right)}(G)$, $S$ has to contain a non-empty zero-sum subsequence. 
\end{proof}

\medskip
Now, in order to obtain effective inequalities from the symbolic constraints of Lemma \ref{lemme clef}, one can use a result proved in \cite{Girard08}, which states that for any finite Abelian group $G$ and every $1 \leq d' \text{ } | \text{ } d \text{ } | \text{ } \exp(G)$, the invariant $\mathsf{D}_{(d',d)}(G)$ is linked with the classical Davenport constant of a particular subgroup of $G$, which can be characterized explicitly. In order to define properly this particular subgroup, we have to introduce the following notation. 

\medskip
For all $i \in \llbracket 1,r \rrbracket$, we set:
$$A_i=\gcd(d',n_i), \text{ }
B_i=\frac{\mathrm{lcm}(d,n_i)}{\mathrm{lcm}(d',n_i)}$$

$$\text{ and } \text{} \upsilon_i(d',d)=\frac{A_i}{\gcd(A_i,B_i)}.$$
For instance, whenever $d$ divides $n_i$, we have $\upsilon_i(d',d)=\gcd(d',n_i)=d'$, and in particular $\upsilon_r(d',d)=d'.$ We can now state our result on $\mathsf{D}_{(d',d)}(G)$ (see \cite{Girard08}, Proposition $3.1$).

\begin{proposition} \label{propmarrantegenerale} Let $G \simeq C_{n_1} \oplus \dots \oplus C_{n_r}$, with $1 < n_1
\text{ } |\text{ } \dots \text{ }|\text{ } n_r \in \mathbb{N}$, be a
finite Abelian group and $d',d \in \mathbb{N}$ be such that $1 \leq d' \text{ } | \text{ } d \text{ } | \text{ } \exp(G)$. Then, we have the following equality:
$$\mathsf{D}_{(d',d)}(G)=\mathsf{D}\left(C_{\upsilon_1(d',d)} \oplus
\dots \oplus C_{\upsilon_r(d',d)}\right).$$
\end{proposition}

\section{Two lemmas related to zero-freeness in $G \simeq C_m \oplus C_{mn}$}
\label{Section lemmes}
In this section, we show how the method presented in Section \ref{Outline of the method} can be used in order to obtain two key lemmas for the proofs of Theorems \ref{main theorem} and \ref{main theorem bis}. To start with, we prove the following result.

\begin{lem}\label{Lemme technique} Let $G \simeq C_m \oplus C_{mn}$, where $m,n \in \mathbb{N}^*$, $n \geq 2$, be a finite Abelian group of rank two, and let $S$ be a zero-sumfree sequence in $G$ with $\left| S \right| \geq \mathsf{d}^*(G)= m+mn-2$. Then, for every $\ell \in \mathcal{D}_n \backslash\{n\}$, one has the following inequality:
$$\displaystyle\sum_{d \in \mathcal{D}_\ell} \alpha_{md} \leq m-1.$$
\end{lem} 

\begin{proof}
Let $S$ be a zero-sumfree sequence in $G \simeq C_m \oplus C_{mn}$ with $\left| S \right| \geq \mathsf{d}^*(G)= m+mn-2$. Let $\ell \in \mathcal{D}_n \backslash\{n\}$, $d'=n/\ell$ and $d=mn$, which leads to $d/d'=m\ell$. We also set $m'=\gcd(d',m)$. Now, let $T$ and $U$ be the two subsequences of $S$ which are defined in Lemma \ref{lemme clef}. In particular, one has $T \subseteq U=S$, and by Proposition \ref{ordreinseqmini} $(i)$, we obtain:
$$\left|T\right|=\displaystyle\sum_{\bar{d} \in \mathcal{D}_\ell} \alpha_{m\bar{d}}.$$

To start with, we determine the exact value of $\mathsf{D}_{(d',d)}(G)$. One has:
\begin{eqnarray*}
\upsilon_1(d',d) & = & \frac{m'}{\gcd\left(m',\frac{\text{lcm}(d,m)}{\text{lcm}(d',m)}\right)}\\
                 & = & \frac{m'}{\gcd\left(m',\frac{d}{d'}\frac{m'}{m}\right)}\\
                 & = & \frac{m'}{\gcd\left(m',m'\ell\right)}\\
                 & = & 1,
\end{eqnarray*}
and, since $\upsilon_2(d',d)=d'$, one obtains, using Proposition \ref{propmarrantegenerale} and Theorem \ref{Proprietes D et Cross number} $(ii)$, the following equalities:

\begin{eqnarray*}
\mathsf{D}_{(d',d)}(G) & = & \mathsf{D}\left(C_{\upsilon_1(d',d)} \oplus C_{\upsilon_2(d',d)}\right)\\
                       & = & \mathsf{D}\left(C_{\frac{n}{\ell}}\right)\\
                       & = & \frac{n}{\ell}.
\end{eqnarray*}

\medskip
Now, let us suppose that one has $\left|T\right| \geq m$. Since $\ell \in \mathcal{D}_n \backslash\{n\}$, we obtain the following inequalities: 
\begin{eqnarray*}
\left|T\right|+\frac{\left|U\right|-\left|T\right|}{\mathsf{D}_{(d',d)}(G)} & \geq & \left|T\right|+\frac{\ell\left(m+mn-2-\left|T\right|\right)}{n}\\ 
& \geq & m+\frac{\ell\left(mn-2\right)}{n}\\
& = & \left(m+m\ell-1\right)-\frac{\ell}{n}+\left(\frac{n-\ell}{n}\right)\\
& > & \left(m+m\ell-1\right)-\frac{\ell}{n}\\
&  = & \mathsf{D}_{\left(\frac{d}{d'},\frac{d}{d'}\right)}(G)-\frac{1}{\mathsf{D}_{(d',d)}(G)},
\end{eqnarray*}
and, according to Lemma \ref{lemme clef}, $S$ must contain a non-empty zero-sum subsequence, which is a contradiction.
Thus, one has $\left|T\right| \leq m-1,$ which is the desired result.
\end{proof} 

Now, let $n \geq 2$ be an integer, and $p_1,\dots,p_r$ be its distinct prime divisors. Given $m \in \mathbb{N}^*$ and a zero-sumfree sequence $S$ in $G \simeq C_m \oplus C_{mn}$ with $\left|S\right|\geq \mathsf{d}^*(G)=m+mn-2$, Lemma \ref{Lemme technique} implies that the integers $\alpha_{md} \in \mathbb{N}$, where $d \in \mathcal{D}_{n}\backslash \{n\}$, have to satisfy the following $r$ linear constraints: $$\displaystyle\sum_{d \in \mathcal{D}_{n/{p_i}}} \alpha_{md} \leq m-1, \text{ for all } i \in \llbracket 1,r \rrbracket.$$  
In the next lemma, we solve a linear integer program on the divisor lattice of $n$, in order to obtain the maximum value of the function 
$$(\alpha_{md})_{d \in \mathcal{D}_{n}\backslash \{n\}} \mapsto \displaystyle\sum_{d \in \mathcal{D}_n \backslash\{n\}}\frac{\alpha_{md}}{d}$$ 
under the $r$ above constraints (the reader interested by linear programming methods is referred to the book \cite{Schri98}, for an exhaustive presentation of the subject).
  
\begin{lem}\label{Lemme technique bis} Let $m,n \in \mathbb{N}^*$, with $n \geq 2$, and let $(x_d)_{d \in  \mathcal{D}_n \backslash\{n\}}$ be a sequence of positive integers, such that for every prime divisor $p$ of $n$, one has the following linear constraint:
$$\displaystyle\sum_{d \in \mathcal{D}_{n/p}} x_d \leq m-1.$$
 Then, one has the following inequality, which is best possible:
$$\displaystyle\sum_{d \in \mathcal{D}_n \backslash\{n\}} \frac{x_{d}}{d} \leq m-1.$$
\end{lem}

\begin{proof} Let $n \geq 2$ be an integer, and $p_1,\dots,p_r$ be its distinct prime divisors. For every $k \in \llbracket 0,m-1 \rrbracket$, let also $\mathcal{S}_k$ be the set of all the sequences of positive integers $x=(x_d)_{d \in  \mathcal{D}_n \backslash\{n\}}$ which verify the above linear constraints, and being such that $x_1=m-k-1$. Now, we can prove, by induction on $k \in \llbracket 0,m-1 \rrbracket$, that the following statement holds. 
$$\text{For every sequence } x \in \mathcal{S}_k, \text{ one has } \displaystyle\sum_{d \in \mathcal{D}_n \backslash\{n\}} \frac{x_{d}}{d} \leq m-1.$$

If $k=0$, then for every $x \in \mathcal{S}_0$, the linear constraints imply that $x_{d}=0$ for all $d \in \mathcal{D}_n \backslash \{1,n\}$, which gives the following equality:
$$\displaystyle\sum_{d \in \mathcal{D}_n \backslash\{n\}} \frac{x_{d}}{d} = m-1.$$ 
Assume now that the statement is valid for $k-1 \geq 0$. Let us define the following map:
\begin{eqnarray*}
f: \mathcal{D}_n \backslash \{n\} & \rightarrow & \{\mathcal{A} \text{ } | \text{ } \emptyset \varsubsetneq \mathcal{A} \subseteq \llbracket 1,r \rrbracket\}\\
d & \mapsto & \{i \in \llbracket 1,r \rrbracket \text{ } | \text{ } d \in \mathcal{D}_{n/p_i}\}.
\end{eqnarray*}
Let $x \in \mathcal{S}_k$ and let $\mathcal{L}$ be the set of the elements $d \in \mathcal{D}_n \backslash \{1,n\}$ such that one has $x_{d} \geq 1$. By definition, and for every $d \in \mathcal{D}_n \backslash \{n\}$, $\left|f(d)\right|$ is the number of linear constraints in which the variable $x_{d}$ appears.  Thus, for every prime divisor $p$ of $n$, $x_{n/p}$ appears in only one linear constraint, and we may assume, without loss of generality, that we have: 
$$\displaystyle\sum_{d \in \mathcal{D}_{n/p}} x_d = m-1.$$ 
Hence, for every $i \in \llbracket 1,r \rrbracket$, the set $\mathcal{L} \cap \mathcal{D}_{n/p_i}$ is non-empty, and one obtains:
$$\displaystyle\bigcup_{d \in \mathcal{L}} f(d) = \llbracket 1,r \rrbracket.$$
Let us consider a non-empty subset $\mathcal{L}'$ of $\mathcal{L}$ verifying the following equality:
$$\displaystyle\bigcup_{d \in \mathcal{L}'} f(d) = \llbracket 1,r \rrbracket,$$
and being of minimal cardinality regarding this property. Since $f(d)$ is a non-empty set for every $d \in \mathcal{D}_n \backslash \{n\}$, the following property has to hold:
$$\displaystyle\bigcup_{d \in \mathcal{L}''} f(d) \subsetneq \llbracket 1,r \rrbracket \text{ for all } \emptyset \subsetneq \mathcal{L}'' \subsetneq \mathcal{L}'.$$
Now, one can notice the following two facts.

\medskip
\noindent 
$\textbf{Fact 1.}$ For every $d \in \mathcal{L}'$, one has $\left|f(d)\right|\leq r-\left|\mathcal{L}'\right|+1$, and in particular, $\left|\mathcal{L}'\right| \leq r$. This fact is a straightforward consequence of the following combinatorial lemma. 

\begin{lem}\label{lemme combinatoire} Let $r \in \mathbb{N}^*$ and $\mathcal{A}_1,\cdots,\mathcal{A}_s$ be $s$ non-empty subsets of $\llbracket 1,r \rrbracket$ verifying
$$\displaystyle\bigcup_{i \in \llbracket 1,s \rrbracket} \mathcal{A}_i = \llbracket 1,r \rrbracket, \text{ and } \displaystyle\bigcup_{i \in I} \mathcal{A}_i \subsetneq \llbracket 1,r \rrbracket \text{ for every subset } \emptyset \subsetneq I \subsetneq \llbracket 1,s \rrbracket.$$
Then, for all $i \in \llbracket 1,s \rrbracket$, one has the following inequality:
$$\left|\mathcal{A}_i\right| \leq r-s+1.$$ 
\end{lem}

\begin{proof} By symmetry, it suffices to prove that one has $\left|\mathcal{A}_1\right| \leq r-s+1$. Assume to the contrary that $\left|\mathcal{A}_1\right| \geq r-s+2$. Since, for all $i \in \llbracket 1,s-1 \rrbracket$, the set $\mathcal{A}_{i+1}$ must contain at least one element from $\llbracket 1,r \rrbracket \backslash \left(\mathcal{A}_1 \cup \cdots \cup \mathcal{A}_i\right)$, one obtains the following inequality: 
$$\left|\mathcal{A}_1 \cup \cdots \cup \mathcal{A}_{i+1}\right| \geq \left|\mathcal{A}_1 \cup \cdots \cup \mathcal{A}_i\right|+1.$$
Therefore, we deduce by an easy induction argument that one has $\left|\mathcal{A}_1 \cup \cdots \cup \mathcal{A}_{s-1}\right| \geq (r-s+2)+(s-2)=r$, and so $\mathcal{A}_1 \cup \cdots \cup \mathcal{A}_{s-1}=\llbracket 1,r \rrbracket$, which is a contradiction. 
\end{proof}

\medskip
\noindent 
$\textbf{Fact 2.}$ For every $d \in \mathcal{D}_n \backslash \{n\}$, one has the following inequalities:
\begin{eqnarray*}
d & \geq & \min f^{-1}(f(d)) \\
     & \geq & \displaystyle\prod_{i \in \llbracket 1,r \rrbracket \backslash f(d)} p^{\nu_{p_i}(n)}_i\\
     & \geq & 2^{r-\left|f(d)\right|}.
\end{eqnarray*} 

Now, using Facts $1$ and $2$, we can prove the desired result, by considering the sequence $y=(y_d)_{d \in  \mathcal{D}_n \backslash\{n\}}$ obtained from $x$ in the following way:
$$
y_d=\begin{cases}
x_1+1 \text{ if } d=1,\\
x_d-1 \text{ if } d \in \mathcal{L}',\\
x_d \hspace{0.75cm} \text{ otherwise}.
\end{cases}
$$
It is readily seen that $y \in \mathcal{S}_{k-1}$. Therefore, Facts $1$ and $2$ give the following inequalities:
\begin{eqnarray*}
m-1 & \geq & \displaystyle\sum_{d \in \mathcal{D}_n \backslash\{n\}} \frac{y_{d}}{d}\\ 
& = & \left(\displaystyle\sum_{d \in \mathcal{D}_n \backslash\{n\}} \frac{x_{d}}{d}\right)+\left(1-\displaystyle\sum_{d \in \mathcal{L}'} \frac{1}{d}\right)\\
& \geq & \left(\displaystyle\sum_{d \in \mathcal{D}_n \backslash\{n\}} \frac{x_{d}}{d}\right)+\left(1-\displaystyle\sum_{d \in \mathcal{L}'} \frac{1}{2^{r-\left|f(d)\right|}}\right)\\
& \geq & \left(\displaystyle\sum_{d \in \mathcal{D}_n \backslash\{n\}} \frac{x_{d}}{d}\right)+\left(1-\displaystyle\sum_{d \in \mathcal{L}'} \frac{1}{2^{\left|\mathcal{L}'\right|-1}}\right)\\
& \geq & \left(\displaystyle\sum_{d \in \mathcal{D}_n \backslash\{n\}} \frac{x_{d}}{d}\right)+\left(1-\frac{\left|\mathcal{L}'\right|}{2^{\left|\mathcal{L}'\right|-1}}\right)\\
& \geq & \displaystyle\sum_{d \in \mathcal{D}_n \backslash\{n\}} \frac{x_{d}}{d},
\end{eqnarray*}
which completes the proof.
\end{proof}

\section{Proofs of the two main theorems}
\label{Demonstrations}
To start with, we show that every finite Abelian group of rank two satisfies Conjecture \ref{Conjecture Cross number}. The following proof of Theorem \ref{main theorem} consists in a direct application of Lemmas \ref{Lemme technique} and \ref{Lemme technique bis}.
\begin{proof}[Proof of Theorem \ref{main theorem}] Let $G \simeq C_m \oplus C_{mn}$, where $m,n \in \mathbb{N}^*$, be a finite Abelian group of rank two, and let $S$ be a zero-sumfree sequence in $G$ with $\left|S\right| \geq \mathsf{d}^*(G)=m+mn-2$. Since, by Theorem \ref{Proprietes D et Cross number} $(ii)$, one has $\mathsf{d}(G)=\mathsf{d}^*(G)$, we obtain that $\left|S\right|=\mathsf{d}^*(G)=m+mn-2$. 

\medskip
If $n=1$, then the desired result follows directly from Proposition \ref{ordreinseqmini} $(i)$, since every element of $S$ has order $m$. Now, let us suppose that $n \geq 2$. Using Proposition \ref{ordreinseqmini} $(i)$, we obtain:
\begin{eqnarray*}
\mathsf{k}(S) & = & \displaystyle\sum_{d \in \mathcal{D}_{mn}} \frac{\alpha_d}{d}\\
              & = & \displaystyle\sum_{d \in \mathcal{D}_n} \frac{\alpha_{md}}{md},
\end{eqnarray*}
and we can distinguish two cases.

\medskip 
$\textbf{Case 1.}$ $\alpha_{mn} \geq mn-1$. In this case, applying Proposition \ref{ordreinseqmini} $(i)$, one obtains:
$$\displaystyle\sum_{d \in \mathcal{D}_n \backslash\{n\}} \alpha_{md}=\left|S\right|-\alpha_{mn},$$ 
which implies the following inequalities:
\begin{eqnarray*}
\mathsf{k}(S) & = &  \displaystyle\sum_{d \in \mathcal{D}_n \backslash\{n\}} \frac{\alpha_{md}}{md} + \frac{\alpha_{mn}}{mn}\\
              & \leq &  \left(\frac{\left|S\right|-\alpha_{mn}}{m}\right) + \frac{\alpha_{mn}}{mn}\\
              & \leq &  \left(\frac{\left|S\right|-(mn-1)}{m}\right) + \left(\frac{mn-1}{mn}\right)\\
              &   =  &  \left(\frac{m-1}{m}\right) + \left(\frac{mn-1}{mn}\right).
\end{eqnarray*} 

\medskip 
$\textbf{Case 2.}$ $\alpha_{mn} \leq mn-1$. Then, by Lemmas \ref{Lemme technique} and \ref{Lemme technique bis}, we obtain:
\begin{eqnarray*}
\mathsf{k}(S) & = &  \displaystyle\sum_{d \in \mathcal{D}_n \backslash\{n\}} \frac{\alpha_{md}}{md} + \frac{\alpha_{mn}}{mn}\\
              & \leq &  \left(\frac{m-1}{m}\right) + \frac{\alpha_{mn}}{mn}\\
              & \leq &  \left(\frac{m-1}{m}\right) + \left(\frac{mn-1}{mn}\right),
\end{eqnarray*} 
which completes the proof.
\end{proof}

Now, we prove Theorem \ref{main theorem bis}, which gives a lower bound for the number of elements with maximal order in a long zero-sumfree sequence of a finite Abelian group of rank two.

\begin{proof}[Proof of Theorem \ref{main theorem bis}] Let $G \simeq C_m \oplus C_{mn}$, where $m,n \in \mathbb{N}^*$, be a finite Abelian group of rank two, and let $S$ be a zero-sumfree sequence in $G$ with $\left|S\right|=\mathsf{d}(G)=m+mn-2$.
\begin{itemize}
\item[$(i)$] Let $p \in \mathcal{P}$ and $a \in \mathbb{N}$ be such that $n=p^a$. If $a=0$, then $G \simeq C_m \oplus C_m$ and, by Proposition \ref{ordreinseqmini} $(i)$, every element of $S$ has order $m$. Now, let us suppose that $a \geq 1$, Then, by Lemma \ref{Lemme technique}, one has:
\begin{eqnarray*}
\left|S\right| - \alpha_{mn} & = & \displaystyle\sum_{d \in \mathcal{D}_{p^{a-1}}} \alpha_{md}\\ 
& \leq & m-1, 
\end{eqnarray*}
which indeed implies that  
\begin{eqnarray*}
\alpha_{mn} & \geq & \left|S\right| - (m-1)\\
            &   =  & m+mn-2-(m-1)\\
            &   =  & mn-1. 
\end{eqnarray*}

\medskip
\item[$(ii)$] 
If $\tau(n) \leq 3$, then $n$ has to be a prime power, and the desired result follows by $(i)$. Now, let us suppose that $\mathcal{D}_n=\left\{d_0=1 < d_1 < d_2 < d_3 \dots\right\}$ contains at least four elements. In particular, one has $n \geq 6$.

\medskip
By Lemmas \ref{Lemme technique} and \ref{Lemme technique bis}, one has:
$$\displaystyle\sum_{d \in \mathcal{D}_n \backslash\{n\}} \frac{\alpha_{md}}{md} \leq \frac{m-1}{m},$$
that is
\begin{eqnarray}\label{relation}
d_1\left(\frac{\alpha_{m\frac{n}{d_1}}}{mn}\right) + d_2\left(\frac{\alpha_{m\frac{n}{d_2}}}{mn}\right) + d_3\left(\frac{\alpha_{m\frac{n}{d_3}}}{mn}\right) + \displaystyle\sum_{\begin{subarray}{c}
d \in \mathcal{D}_n\\
d > d_3\end{subarray}} d\left(\frac{\alpha_{m\frac{n}{d}}}{mn}\right) \leq \frac{m-1}{m}.
\end{eqnarray}
Now, we can distinguish two cases.

\medskip
\textbf{Case 1.} $d_3=4$. Then $d_1=2$, $d_2=3$ and (\ref{relation}) implies:
\begin{eqnarray}\label{relation bis}
2\left(\frac{\alpha_{m\frac{n}{2}}}{mn}\right) + 3\left(\frac{\alpha_{m\frac{n}{3}}}{mn}\right) + 4\left(\frac{\alpha_{m\frac{n}{4}}}{mn}\right) + 5\left(\displaystyle\sum_{\begin{subarray}{c}
d \in \mathcal{D}_n\\
d > d_3\end{subarray}} \frac{\alpha_{m\frac{n}{d}}}{mn}\right) \leq \frac{m-1}{m}.
\end{eqnarray}
But since 
$$\displaystyle\sum_{\begin{subarray}{c}
d \in \mathcal{D}_n\\
d > d_3\end{subarray}} \alpha_{m\frac{n}{d}}=\left|S\right|-\alpha_{mn}-\alpha_{m\frac{n}{2}}-\alpha_{m\frac{n}{3}}-\alpha_{m\frac{n}{4}},$$
relation (\ref{relation bis}) implies:
$$5\left(\frac{m+mn-2-\alpha_{mn}}{mn}\right)-3\left(\frac{\alpha_{m\frac{n}{2}}}{mn}\right) -2\left(\frac{\alpha_{m\frac{n}{3}}}{mn}\right) -\left(\frac{\alpha_{m\frac{n}{4}}}{mn}\right) \leq \frac{m-1}{m},$$
that is
$$5\left(\frac{m+mn-2-\alpha_{mn}}{mn}\right)-\left(\frac{\alpha_{m\frac{n}{2}}+\alpha_{m\frac{n}{4}}}{mn}\right) -2\left(\frac{\alpha_{m\frac{n}{2}}+\alpha_{m\frac{n}{3}}}{mn}\right)\leq \frac{m-1}{m},$$
Now using the fact that, by Lemma \ref{Lemme technique}, one has:
$$\alpha_{m\frac{n}{2}}+\alpha_{m\frac{n}{4}} \leq m-1 \text{ } \text{ as well as } \text{ } \alpha_{m\frac{n}{3}} \leq m-1,$$
we obtain
$$5\left(\frac{m+mn-2-\alpha_{mn}}{mn}\right)-\left(\frac{m-1}{mn}\right) -2\left(\frac{2\left(m-1\right)}{mn}\right) \leq \frac{m-1}{m},$$
which is equivalent to
$$5\left(\frac{mn-1-\alpha_{mn}}{mn}\right) \leq \frac{m-1}{m},$$
that is
$$5(mn-1) - n(m-1) \leq 5\alpha_{mn},$$
and thus
$$\frac{4}{5}mn+\frac{(n-5)}{5} \leq \alpha_{mn},$$
which is the desired result.

\medskip
\textbf{Case 2.} $d_3 \geq 5$. Then (\ref{relation}) implies:
\begin{eqnarray}\label{relation ter}
d_1\left(\frac{\alpha_{m\frac{n}{d_1}}}{mn}\right) + d_2\left(\frac{\alpha_{m\frac{n}{d_2}}}{mn}\right) + 5\left(\displaystyle\sum_{\begin{subarray}{c}
d \in \mathcal{D}_n\\
d \geq d_3\end{subarray}} \frac{\alpha_{m\frac{n}{d}}}{mn}\right) \leq \frac{m-1}{m}.
\end{eqnarray}
But since 
$$\displaystyle\sum_{\begin{subarray}{c}
d \in \mathcal{D}_n\\
d \geq d_3\end{subarray}} \alpha_{m\frac{n}{d}}=\left|S\right|-\alpha_{mn}-\alpha_{m\frac{n}{d_1}}-\alpha_{m\frac{n}{d_2}},$$
relation (\ref{relation ter}) implies:
$$5\left(\frac{m+mn-2-\alpha_{mn}}{mn}\right)+(d_1-5)\left(\frac{\alpha_{m\frac{n}{2}}}{mn}\right) + \left(d_2-5\right)\left(\frac{\alpha_{m\frac{n}{d_2}}}{mn}\right) \leq \frac{m-1}{m}.$$
Therefore, since $d_1 \geq 2$ and $d_2 \geq 3$, we have
$$5\left(\frac{m+mn-2-\alpha_{mn}}{mn}\right)-3\left(\frac{m-1}{mn}\right) - 2\left(\frac{m-1}{mn}\right) \leq \frac{m-1}{m},$$
that is
$$5\left(\frac{mn-1-\alpha_{mn}}{mn}\right) \leq \frac{m-1}{m},$$
which leads to 
$$\frac{4}{5}mn+\frac{(n-5)}{5} \leq \alpha_{mn},$$
and the proof is complete.
\end{itemize}
\end{proof}

\section{A concluding remark}
\label{concluding remark}
Given a finite Abelian group $G$, the investigation of the maximal possible length of a zero-sumfree sequence $S$ in $G$ with large cross number may also be of interest. Concerning this question, we propose the following general conjecture, which can be seen as a dual version of Conjecture \ref{Conjecture Cross number}.

\begin{conjecture}\label{Conjecture Davenport} 
Let $G$ be a finite Abelian group and $G \simeq C_{\nu_1} \oplus \dots \oplus C_{\nu_s}$, with $\nu_i > 1$ for all $i \in \llbracket 1,s \rrbracket$, be its longest possible decomposition into a direct product of cyclic groups. 
Given a zero-sumfree sequence $S$ in $G$ verifying $\mathsf{k}(S) \geq \mathsf{k}^*(G)$, one always has the following inequality:
\begin{eqnarray*} \left| S \right| \leq \displaystyle\sum_{i=1}^s (\nu_i-1). \end{eqnarray*}
\end{conjecture}

It can easily be seen, by Theorem \ref{Proprietes D et Cross number} $(i)$, that Conjecture \ref{Conjecture Davenport} holds true for finite Abelian $p$-groups. Even in the case of finite cyclic groups which are not $p$-groups, this problem is still wide open. Yet, in this special case, the following result supports the idea that a zero-sumfree sequence with large cross number has to be a "short" sequence.

\begin{theorem} \label{theorem index} Let $n \in \mathbb{N}^*$ be such that $n$ is not a prime power, and let $S$ be a zero-sumfree sequence in $C_n$ verifying $\mathsf{k}(S) \geq \mathsf{k}^*(C_n)$. Then, one has the following inequality: 
$$\left|S\right| \leq \left\lfloor \frac{n}{2}\right\rfloor.$$ 
\end{theorem}

\begin{proof} So as to prove this result, we will use the notion of \em index \em of a sequence in a finite cyclic group, which was introduced implicitly in \cite{KleitLemke89}, Conjecture p.$344$, and more explicitly in \cite{ChapFreezeSmith99}. Let $g \in C_n$ with $\text{ord}(g)=n$, and let $S=(g_1,\dots,g_{\ell})=(n_1g,\dots,n_{\ell}g)$, where $n_1,\dots,n_{\ell} \in \llbracket 0,n-1 \rrbracket$, be a sequence in $C_n$. We define:
$$\left\|S\right\|_g=\displaystyle\sum^{\ell}_{i=1}\frac{n_i}{n}.$$  
Since, for every $i \in \llbracket 1,\ell \rrbracket$, we have 
$$\frac{\gcd(n_i,n)}{n}=\frac{1}{\text{ord}(g_i)},$$
one can notice that $\left\|S\right\|_g \geq \mathsf{k}(S)$ for all $g \in C_n$ with $\text{ord}(g)=n$. Then, the index of $S$, denoted by $\text{index}(S)$, is defined in the following fashion:
$$\text{index}(S)=\displaystyle\min_{\tiny{\begin{array}{c}
g \in C_n \\
\text{ord}(g)=n
\end{array}}} \left\|S\right\|_g.$$

Now, if $n$ is not a prime power and $S$ is a zero-sumfree sequence in $C_n$ such that $\mathsf{k}(S) \geq \mathsf{k}^*(C_n)$, one obtains, by the very definition of the index, the following inequalities: 
\begin{eqnarray*} \text{index}(S)    &  \geq  & \mathsf{k}(S) \\
                                     & \geq & \mathsf{k}^*(C_n)\\ 
                                     &   >  &  1.
\end{eqnarray*}
Therefore, using a result of Savchev and Chen (see Theorem $9$ in \cite{SavChen07}), one must have the following inequality: 
$$\left|S\right| \leq \left\lfloor \frac{n}{2}\right\rfloor,$$
which completes the proof.
\end{proof}

In particular, Theorem \ref{theorem index} implies that Conjecture \ref{Conjecture Davenport} holds true for all the cyclic groups of the form $C_{2p^a}$, where $p \in \mathcal{P}$ and $a \in \mathbb{N}$.

\section*{Acknowledgments}
I am grateful to my Ph.D. advisor Alain Plagne for his help during the preparation of this paper. I would like also to thank Alfred Geroldinger and Wolfgang Schmid for useful remarks on a preliminary version of this work.


\end{document}